\documentclass[12pt]{amsart}
\usepackage{amssymb,yhmath,latexsym,amsmath,amscd}
\theoremstyle{plain}
\newtheorem{thm}{Theorem}[section]
\newtheorem{prop}[thm]{Proposition}
\newtheorem{cor}[thm]{Corollary}
\newtheorem{lemma}[thm]{Lemma}

\theoremstyle{definition}
\newtheorem{defn}[thm]{Definition}
\newtheorem{rmk}[thm]{Remark}

\newtheorem{notat}[thm]{Notation}

\newtheorem{ex}[thm]{Example}

\newcommand{\lra}{\longrightarrow}
\newcommand{\PP}{\mathbb{P}}

\newcommand{\ZZ}{\mathbb{Z}}

\newcommand{\OO}{\mathcal{O}}

\newcommand{\K}{\mathcal{K}}
\newcommand{\I}{\mathcal{I}}
\newcommand{\M}{\mathcal{M}}

\newcommand{\MM}{\mathfrak{m}}

\newcommand{\FF}{\mathbb{F}}
\newcommand{\GG}{\mathbb{G}}

\newcommand{\OP}[1]{\OO_{\PP^{#1}}}

\newcommand{\Lra}{\Longrightarrow}
\newcommand{\hgt}{\mbox{ht}\;}

\title{The G-biliaison class of symmetric determinantal schemes}

\author{Elisa Gorla}

\address{Institut f\"ur Mathematik
\\ Universit\"at Z\"urich, \hfil\break\indent Winterthurerstrasse
190, CH-8057 Z\"urich, Switzerland}

\email{elisa.gorla@math.unizh.ch}

\begin{document}

\maketitle

{\bf Abstract:} We consider a family of schemes, that are defined by minors 
of a homogeneous symmetric matrix with polynomial entries.
We assume that they have 
maximal possible codimension, given the size of the
matrix and of the minors that define them. 
We show that these schemes are G-bilinked to a
linear variety of the same dimension. In particular, they can be 
obtained from a linear
variety by a finite sequence of ascending G-biliaisons on some
determinantal schemes. 
We describe the biliaisons explicitly in the proof of
Theorem~\ref{bilsym}. 
In particular, it follows that 
these schemes are glicci.

\begin{center}{\bf Introduction}\end{center}

A main open question in liaison theory consists of deciding 
whether every arithmetically Cohen-Macaulay scheme is in the 
same G-liaison class of a complete intersection (glicci). 
The Theorem of Gaeta says that every arithmetically
Cohen-Macaulay scheme of codimension 2 belongs to the CI-liaison class of
a complete intersection. This result was generalized by
Kleppe, Migliore, Mir\'o-Roig, Nagel, and Peterson in
\cite{kl01}, where they proved that every standard determinantal
scheme is in the (even) G-liaison class of a complete
intersection of the same codimension. Hartshorne strengthened
their result in \cite{ha04u2}, proving that every standard
determinantal scheme is in the G-biliaison class of a complete intersection. 
In \cite{ca01a}, Casanellas and
Mir\'o-Roig proved that any arithmetically Cohen-Macaulay divisor
on a rational normal scroll surface is glicci. In her
Ph.D. thesis \cite{ca02t},
M.~Casanellas generalized this result to arithmetically
Cohen-Macaulay divisors on a rational normal scroll (see also \cite{ca01a2} 
for a statement of the same result).
It easily follows from the main
theorem in the paper of Watanabe \cite{wa73a} that every arithmetically 
Gorenstein scheme of codimension 3 is licci (i.e. it belongs to the CI-liaison 
class of a complete intersection). Moreover, M.~Casanellas, E.~Drozd and 
R.~Hartshorne recently proved that every arithmetically Gorenstein scheme 
is glicci, regardless of its codimension (see Theorem~7.1 of~\cite{ca04u}). 
We are still far from being able to answer in full generality to the question 
of whether any arithmetically Cohen-Macaulay scheme is glicci. 
As far as we know, all the work in this direction deals with
specific families of arithmetically-Cohen Macaulay schemes.
In this note we consider a family of schemes whose saturated
ideal is generated by minors of a fixed size of a symmetric
matrix with polynomial entries. We prove that these schemes are
G-bilinked to a complete intersection. In particular, they are glicci.

In the first section, we introduce the family of schemes that
will be the object of our study. 
Their defining ideals are generated by the minors of a fixed size 
of a symmetric matrix with polynomial entries. We assume that the 
schemes have the
highest possible codimension, for a fixed size of the matrix and
of the minors that define them.
We call a scheme of this kind {\bf symmetric determinantal}. 
We observe that symmetric 
determinantal schemes do not exist for any given codimension, 
in fact they can only
have codimension $b\choose 2$, for some $b\geq 2$. Examples
of symmetric determinantal schemes are complete intersections of 
admissible codimension, the Veronese surface in $\PP^5$, and some
standard and good determinantal schemes (see Examples~\ref{ci},
\ref{vero}, and~\ref{indet}). In the first section we also
introduce the concept of {\bf almost-symmetric determinantal
  scheme}. Symmetric and
almost-symmetric determinantal schemes are arithmetically
Cohen-Macaulay. This follows from a result of Kutz \cite{ku74a}.
In Theorem~\ref{subm} we discuss when an
almost-symmetric determinantal scheme is
generically complete intersection, by giving an equivalent
condition and a sufficient condition. In Theorem~\ref{subsd} we
state the analogous result for symmetric determinantal schemes.
In Theorem~\ref{ht1} we give an upper bound on the
height of the ideal of minors of size $t\times t$ of a symmetric
$m\times m$ matrix modulo the ideal of minors of size $t\times t$
of the same matrix that do not involve the last row.

Proposition~\ref{gsd} in Section 2 clarifies the
connection between symmetric and almost-symmetric determinantal
schemes. For each symmetric determinantal scheme $X$ we produce
an almost-symmetric determinantal scheme $Y$ such that $X$ is a
generalized divisor on $Y$, $Y$ is arithmetically Cohen-Macaulay and 
generically complete intersection. We also
construct another symmetric determinantal scheme $X'$
that is a generalized divisor on $Y$. 
Theorem~\ref{bilsym} is the main result of this paper:
any symmetric determinantal scheme belongs to the
G-biliaison class of a linear variety. All the divisors involved
in the G-biliaisons are symmetric determinantal, and the
G-biliaisons are performed on almost-symmetric determinantal
schemes.

We wish to emphasize the analogy from the point 
of view of liaison theory between the family of 
symmetric determinantal schemes and the family 
of standard determinantal schemes. See \cite{kr00} for 
the definition of standard and good determinantal
schemes. See also \cite{kl01} for some of their
properties, mainly in relation with liaison theory.
In \cite{kl01} Kleppe, Migliore, Mir\`o-Roig, Nagel and Peterson 
prove that standard determinantal schemes are glicci. Their argument
is constructive, meaning that following the proof of Theorem 3.6 
of \cite{kl01} one can write down explicitly all the links.  
In \cite{ha04u2}, Hartshorne proves that standard determinantal 
schemes are in the same G-biliaison class of a linear variety.
He shows that the G-bilinks constructed in \cite{kl01} can indeed be 
regarded as elementary G-biliaisons.
Our main theorem is analogous to the main theorem in \cite{ha04u2}, and 
the G-bilinks that one obtains following our proof are in the 
spirit of \cite{kl01}.

After submission of this paper, the author was informed by Kleppe and Mir\'o-Roig 
that they independently proved that the ideal generated by the 
submaximal minors of a symmetric square matrix is glicci. They work under the 
assumption that the codimension of the ideal is 3 (hence maximal). The result is 
unpublished, and it is part of the work in progress \cite{kl05u}.
Notice that the ideal of submaximal minors of a square matrix that is not symmetric 
is Gorenstein (hence glicci) whenever it has maximal codimension 4.
 
{\bf Acknowledgment:} I wish to express my gratitude 
to Juan Migliore and Claudia Polini for their kind invitation and
hospitality during my visit at the University of Notre Dame.
I also wish to thank them for useful comments and discussions about the 
material of this paper. I am grateful to Aldo Conca for useful
comments on the first draft of this paper, and for finding a better
proof of Lemma~\ref{det&ci}.

\section{Symmetric and almost-symmetric determinantal schemes}

Let $X$ be a scheme in $\PP^n=\PP^n_K$, where $K$ is an algebraically
closed field. We assume that the characteristic of $K$ is different
from $2$. 
Let $I_X$ be the saturated homogeneous ideal corresponding to $X$ in the 
polynomial ring $R=K[x_0,x_1,\ldots,x_n]$. We denote by $\MM$ the
homogeneous irrelevant maximal ideal of $R$, $\MM=(x_0, x_1, \ldots ,x_n)$.
For an ideal $I\subset R$, we denote by $H^0_*(I)$ the saturation of $I$ 
with respect to the maximal ideal $\MM$.

Let $\I_X\subset \OP{n}$ be the ideal sheaf of~$X$.
Let $Y$ be a scheme that contains $X$. We denote by $\I_{X|Y}$ the ideal sheaf 
of $X$ restricted to $Y$, i.e. the quotient sheaf $\I_X/\I_Y$.
For $i\geq 0$, we let
$H^i_*(\PP^n,\I)=\oplus_{t\in\ZZ}H^i(\PP^n,\I(t))$ denote the
$i$-th cohomology module
of the sheaf $\I$ on $\PP^n$. We write just $H^0_*(\I)$, when the 
ambient space $\PP^n$ is clearly defined.

\begin{notat} Let $I\subseteq R$ be a homogeneous ideal. We let
  $\mu(I)$ denote the cardinality of a set of minimal generators
  of $I$.
\end{notat}

In this paper we deal with schemes whose saturated ideals are 
generated by minors of matrices with polynomial entries. We refer the reader 
to \cite{kr00} for the definition of standard and good determinantal
schemes.

\begin{defn}
Let $M$ be a matrix of size $m\times m$ with entries in $R$. We say that $M$ 
is {\em t-homogeneous} if the minors of $M$ of size $s\times s$ are 
homogeneous polynomials for all $s\leq t$. We say that $M$ is 
{\em homogeneous} if it is m-homogeneous.
\end{defn}

We will always consider t-homogeneous matrices. Moreover, we will 
regard symmetric matrices up to invertible linear transformations that
preserve their symmetry, and almost symmetric matrices  
up to invertible linear transformations that
preserve the property of being almost-symmetric. We regard all
matrices up to changes of coordinates.
See Definition~\ref{asm} for the definition of almost-symmetric matrix.

\begin{defn}\label{sd}
Let $X\subset\PP^n$ be a scheme. We say that $X$ is 
{\em symmetric determinantal} if: 
\begin{enumerate}
\item there exists a symmetric t-homogeneous matrix
$M$ of size $m\times m$ with entries in $R$, such that the saturated 
ideal of $X$ is
generated by the minors of size $t\times t$ of $M$, $I_X=I_t(M)$,
and
\item $X$ has codimension $m-t+2\choose 2$.
\end{enumerate}
\end{defn}

\begin{rmk}\label{codim}
For any scheme $X$ satisfying requirement (1) of
Definition~\ref{sd}, we have $$codim(X)\leq {m-t+2\choose 2}.$$
See Theorem 2.1 of \cite{jo78} for a proof of this fact.
\end{rmk}

The remark shows that symmetric determinantal schemes have highest 
possible codimension among the schemes defined by minors of a symmetric
matrix, for a given size of the matrix and of the minors.
Notice also that by requiring that $I_t(M)$ is the defining ideal of 
a scheme $X$, we are requiring that the ideal of $t\times t$ minors
of $M$ is saturated.

It is worth emphasizing that symmetric determinantal schemes do not occur 
for every codimension. However, complete intersections are a special case of 
symmetric determinantal schemes for the codimensions for which symmetric 
determinantal schemes do exist.

\begin{rmk}
Symmetric determinantal schemes do not exist for every
codimension. In fact, we have symmetric determinantal
schemes of codimension $c$ if and only if $c$ is of the form $b\choose
2$, for some integer $b\geq 2$.
\end{rmk}

Complete intersections are an easy example of symmetric
determinantal schemes, for each admissible codimension.

\begin{ex}\label{ci}
Let $X$ be a complete intersection of codimension ${b\choose 2}$.
Then
$$I_X=(F_{ij} \; | \; 1\leq i\leq j\leq b-1 )=I_1(M)$$
where $M=(F_{ij}')$, $F_{ij}'=F_{ij}$ if $i\leq j$ and
$F_{ij}'=F_{ji}$ if $i\geq j$.
The matrix $M$ is symmetric of size $(b-1)\times(b-1)$. Its
entries, i.e. its minors of size $1$ define
$X$, and the codimension of $X$ is ${b\choose 2}$.
Hence $X$ is symmetric determinantal.
\end{ex}

The Veronese surface $V\subset\PP^5$ is an example of a symmetric
determinantal scheme that is not a complete intersection.

\begin{ex}\label{vero}
Let $V\subset\PP^5$ be the Veronese surface. The saturated
ideal of $V$ is minimally generated by the (distinct) minors of
size two by two of a symmetric matrix of indeterminates of size 
three by three:
$$I_V=I_2 \left[\begin{array}{ccc} 
x_0 & x_1 & x_2 \\
x_1 & x_5 & x_3 \\
x_2 & x_3 & x_4
\end{array}\right].$$
The Veronese surface has codimension $3={3-2+2\choose 2}$.
\end{ex}

Proposition 2.1 of \cite{ku74a} provides an example of a symmetric
determinantal scheme in $\PP^n$ for $n={m+1\choose 2}$ and for each $t\leq m$. 

\begin{ex}\label{indet}
For any fixed $m\geq 1$, and for any choice of $t$ such that
$1\leq t\leq m$, let $n={m+1\choose 2}$.
Let $X\subset\PP^n$
be the symmetric determinantal scheme whose saturated
ideal is generated by the minors of
size $t\times t$ of the symmetric matrix of indeterminates of size 
$m\times m$:
$$I_X=I_t \left[\begin{array}{cccc} 
x_{1,1} & x_{1,2} & \cdots & x_{1,m} \\
x_{1,2} & x_{2,2} & \cdots & x_{2,m} \\
\vdots & \vdots & & \vdots \\
x_{1,m} & x_{2,m} & \cdots & x_{m,m}
\end{array}\right].$$
From 
From Proposition 2.1 in \cite{ku74a} we have that
$X$ has $$codim(X)=depth(I_X)={m-t+2\choose 2},$$
therefore $X$ is arithmetically Cohen-Macaulay and symmetric determinantal.
\end{ex}

In his dissertation \cite{co93t}, A. Conca studied the ideals of
Example~\ref{indet}, and in fact a larger family of ideals
generated by minors of matrices of indeterminates.
He computed the Gr\"obner basis of the ideals of
Example~\ref{indet} with respect to a diagonal
monomial order. As a consequence, he was able to compute some of
the invariants related to the Poincar\'e series of these
rings. He also showed that the schemes corresponding to these
ideals are reduced, irreducible and normal, and he characterized
the arithmetically Gorenstein ones among them. See section 4 of \cite{co93t}
for more details.

Complete intersection schemes of codimension ${b\choose 2}$ for some $b$ 
are good determinantal schemes that are also symmetric
determinantal (as observed in Example~\ref{ci}).
Notice however that symmetric determinantal schemes are not a subfamily of
standard or good determinantal schemes.
For example, the Veronese surface in $\PP^5$ is a symmetric determinantal 
scheme, but it is not standard determinantal (see Proposition 6.7 in \cite{go04t}).
Moreover, in Proposition 6.7 and Proposition 6.17 of \cite{go04t}
we provide a large class of examples of symmetric determinantal schemes 
that are not standard determinantal. They include the schemes of 
Example~\ref{indet}. Other examples are analyzed in \cite{go04u}.

\begin{defn}\label{asm}
Let $O$ be a matrix of size $(m-1)\times m$. We say that $O$ is 
{\em almost symmetric} if
the submatrix of $O$ consisting of the first $m-1$ columns
is symmetric.
\end{defn}

\begin{defn}\label{asd}
Let $Y\subset\PP^n$ be a scheme. We say that $Y$ is 
{\em almost-symmetric determinantal} if: 
\begin{enumerate}
\item there exists an almost-symmetric t-homogeneous matrix
$O$ of size $(m-1)\times m$ with entries in $R$, such that the 
saturated ideal of $Y$ is
generated by the minors of size $t\times t$ of $O$, $I_Y=I_t(O)$
\item $Y$ has codimension ${m-t+2\choose 2}-1$.
\end{enumerate}
\end{defn}

Notice that in analogy with the case of symmetric determinantal schemes, 
the ideal of $t\times t$ minors of $O$ is saturated since it is the 
ideal associated to a projective scheme.
Moreover, we require almost-symmetric determinantal schemes to have highest 
possible codimension among the schemes defined by minors of an almost-symmetric
matrix, for a given size of the matrix and of the minors.

\begin{rmk}\label{codas}
For any scheme $Y$ satisfying requirement (1) of
Definition~\ref{asd}, we have $$codim(Y)\leq {m-t+2\choose
  2}-1.$$ See for example the paper of Kutz \cite{ku74a}.
A. Conca showed in his Ph.D. dissertation that if $O$ is an
almost-symmetric matrix of indeterminates, then
$$codim(Y)={m-t+2\choose 2}-1$$
(see Proposition~4.6.2 of~\cite{co93t}).
\end{rmk}

The previous remark provides us with an example of
almost-symmetric determinantal schemes.

\begin{ex}\label{as-indet}
Let $O$ be an almost-symmetric matrix of indeterminates of size\newline
$(m-1)\times m$. Let $n={m+1\choose 2}$.
For any choice of $1\leq t\leq m-1$, let $Y_t\subseteq\PP^n$ be 
the scheme whose
saturated ideal is generated by the minors of size $t\times t$ of
$M$. Then it follows from Proposition~4.6.2 of~\cite{co93t} that
$$codim(Y_t)={m-t+2\choose 2}-1.$$
Hence $Y_t$ is an almost-symmetric determinantal scheme.
\end{ex}

Similarly to the case of symmetric determinantal schemes,
almost-symmetric determinantal schemes do not exist for any
codimension. This is clear from part (2) of Definition~\ref{asd}.
Notice also that complete intersections of codimension
${b\choose 2}-1$ for some $b\geq 3$ are almost-symmetric determinantal.

\begin{ex} 
Every complete intersection of codimension ${b\choose 2}-1$ for
some $b\geq 3$ is an almost-symmetric determinantal scheme. 
In fact, let $Y\subset\PP^n$ be a complete intersection of
codimension ${b\choose 2}-1$. Notice that 
${b\choose 2}-1={b-1\choose 2}+b-2$. 
Let $$I_Y=(F_{i,j}, G_k\; |\; 1\leq i\leq j\leq b-2,\; 1\leq
k\leq b-2)$$ be a minimal system of generators of the ideal of
$Y$. Then $I_Y$ is generated by the entries of the 
almost-symmetric matrix $$O=\left[\begin{array}{ccccc}
F_{1,1} & F_{1,2} & \cdots & F_{1,b-2} & G_1 \\
F_{1,2} & F_{2,2} & \cdots & F_{2,b-2} & G_2 \\
\vdots & \vdots & & \vdots & \vdots \\
F_{1,b-2} & F_{2,b-2} & \cdots & F_{b-2,b-2} & G_{b-2} \\
\end{array}\right],$$
hence $Y$ is almost-symmetric determinantal.
\end{ex}

\begin{rmk}
The family of almost-symmetric determinantal schemes does
not coincide with the family of symmetric determinantal schemes,
since by Remark~\ref{codim} and Remark~\ref{codas} it follows
that in general they have different codimensions.
\end{rmk}

Almost-symmetric determinantal schemes are not
a subfamily of standard or good determinantal
schemes. The schemes of Example~\ref{as-indet} 
are a family of almost-symmetric determinantal 
schemes that are not standard determinantal whenever $1<t<m$.

Cohen-Macaulayness of symmetric and almost-symmetric 
determinantal schemes was proved by R. Kutz in Theorem 1 of \cite{ku74a}.
We state here a special case of this result, as we will 
need it in this section.

\begin{thm}[Kutz]\label{gsd1}
Symmetric and almost-symmetric determinantal schemes are 
arithmetically Cohen-Macaulay.
\end{thm}

We now establish some further properties of almost-symmetric 
determinantal schemes that will be needed in this paper. 
We use the notation of Definition~\ref{asd}.

We start by observing that a scheme defined by the $t\times t$
minors of a t-homogeneous matrix
is a complete intersection only when it is generated by the
entries of the matrix, or by its determinant (in the case of a 
square matrix). In \cite{go79a}, Goto has shown that the ideal of $t\times t$
minors of a symmetric matrix of indeterminates of size $m\times
m$ is Gorenstein if and only if $t=1$ or $m-t$ is even. 
However, here we need a result that applies to almost-symmetric
matrices, and in general to matrices with polynomial entries.

\begin{lemma}\label{det&ci}
Let $M$ be a t-homogeneous symmetric matrix of size $m\times m$,
or a t-homogeneous almost-symmetric matrix of size $(m-1)\times
m$, $1\leq t\leq m-1$.
Let $M$ have entries in $R$ or in $R_P$ for some prime $P$, and 
assume that $M$ has no invertible entries. 
If $I_t(M)$ is a complete intersection, then $t=1$.
\end{lemma}

\begin{proof}
We know that the thesis is true for $M=(x_{ij})$ a generic
symmetric or almost-symmetric matrix. In fact, in both cases a
minimal system of generators of $M$ is given by the minors of the
$t\times t$ submatrices of $M$ whose diagonal is on or above the
diagonal of $M$. Comparing the number of such minors and the
codimension of $I_t(M)$, we conclude that they agree if and
only if $t=1$. 

Consider now the general case when $M=(F_{ij})$ has entries in $R$ or
$R_P$ and $I_t(M)$ is a complete intersection. Let
$I=(x_{ij}-F_{ij})_{ij}\subseteq R[x_{ij}]$ and
let $S=R[x_{ij}]/I$. Let $N=(x_{ij})$ be a symmetric or almost-symmetric
matrix of indeterminates of the same size as $M$. If we tensor a
minimal free resolution of $I_t(N)$ by $S$, we obtain a free
resolution for $I_t(M)$.
This follows from Theorem 3.5 in \cite{br88b}.
The resolution that we obtain for $I_t(M)$ is
minimal, since we construct it from a minimal free resolution by
substituting the indeterminates $x_{ij}$ with $F_{ij}$. Therefore,
there can be no invertible entry in any of the maps. Then
$I_t(M)$ is a complete intersection if and only if $I_t(N)$ is a
complete intersection, and the thesis follows.
\end{proof}

\begin{defn}
Let $X\subset\PP^n$ be a scheme. $X$ is {\em generically complete
  intersection} if the localization $(I_X)_P$ is a complete
intersection for every $P$ minimal associated prime of $I_X$.

$X$ is {\em generically Gorenstein}, abbreviated $G_0$, if the
localization $(I_X)_P$ is a Gorenstein ideal
for every $P$ minimal associated prime of $I_X$.
\end{defn}

The next theorem will be used in the proof of
Theorem~\ref{bilsym}.

\begin{thm}\label{subm}
Let $Y$ be an almost-symmetric determinantal scheme with defining
matrix $O$, $I_Y=I_t(O)$. Let $c={m-t+2\choose 2}-1$ 
be the codimension of $Y$. The following are equivalent:
\begin{enumerate}
\item $Y$ is generically complete intersection.
\item $\hgt I_{t-1}(O)\geq c+1$.
\end{enumerate}
Let $N$ be the symmetric matrix obtained from $O$ by
  deleting the last column.
The two equivalent conditions are verified if $\hgt I_{t-1}(N)=c+1$.
\end{thm}

\begin{proof}
\begin{enumerate} 

\item $\Lra$ (2): since $I_t(O)\subseteq I_{t-1}(O)$, then $\hgt
I_{t-1}(O)\geq c$. Therefore it suffices to show that $\hgt
I_{t-1}(O)\neq c$. By contradiction, assume that there exists a
minimal associated prime $P$ of $I_{t-1}(O)$ of height $c$. Then
$P$ is also a minimal associated prime of $I_t(O)$.
Let 
$$\begin{array}{rcl}
\varphi: \FF & \lra & \GG \\
v & \mapsto & Ov \end{array}$$
be the map induced by $O$. Here $v$ is a column vector whose entries
are polynomials, $\FF$ and $\GG$ are free $R$-modules of ranks
$m$ and $m-1$ respectively.
By Proposition~16.3 in~\cite{br88b} we have that the map
$\varphi_P$ that we obtain from $\varphi$ after localizing at the
prime ideal $P$ is an isomorphism on a direct summand $R_P^s$
of $\FF_P$, $\GG_P$ for some $s\leq t-1$. We let $s$ be maximal
with this property.
The localization $O_P$ of $O$ at $P$
can be reduced after invertible row and column operations to
the form $$O_P=\left[\begin{array}{cc}
I_s & 0 \\ 0 & B
\end{array}\right],$$
where $I_s$ is an identity matrix of size $s\times s$, 
$0$ represents a matrix of zeroes,
and $B$ is a matrix of size $(m-s)\times (m-1-s)$ that has no
invertible entries. By assumption, $I_t(O)_P\subseteq R_P$ is a
complete intersection ideal. Since 
$$I_t(O)_P=I_t(O_P)=I_{t-s}(B)$$
and $B$ has no invertible entries, it follows by Lemma~\ref{det&ci}
that $t-s=1$, that is $s=t-1$. But then
$$I_{t-1}(O)_P=I_{t-1}(O_P)=R_P,$$ that contradicts the
assumption that $P\supseteq I_{t-1}(O)$. 

\item $\Lra$ (1): let $P$ be a minimal associated prime of
$I_t(O)$. Then $$\hgt P=c<c+1\leq\hgt I_{t-1}(O)$$ so
$P\not\supseteq I_{t-1}(O)$. Let 
$$\begin{array}{rcl}
\varphi: \FF & \lra & \GG \\
v & \mapsto & Ov \end{array}$$
be the map induced by $O$. Here $v$ is a column vector whose entries
are polynomials, $\FF$ and $\GG$ are free $R$-modules of ranks
$m$ and $m-1$ respectively.
By Proposition~16.3 in~\cite{br88b} we have that the map
$\varphi_P$ that we obtain from $\varphi$ after localizing at the
prime ideal $P$ is an isomorphism on a direct summand $R_P^{t-1}$
of $\FF_P$, $\GG_P$. Hence, the localization $O_P$ of $O$ at $P$
can be reduced, after elementary row and column operations, to
the form $$O_P=\left[\begin{array}{cc}
I_{t-1} & 0 \\ 0 & B
\end{array}\right],$$
where $I_{t-1}$ is an identity matrix of size $(t-1)\times
(t-1)$, $0$ represents a matrix of zeroes,
and $B$ is a matrix of size $(m+1-t)\times (m-t).$ We claim that
$B$ is an almost-symmetric matrix. This is clear if the symmetric
part of $O_P$ contains an invertible minor of size $(t-1)\times
(t-1)$. If instead we have an invertible $(t-1)\times (t-1)$
minor of $O$ that contains the last column, we can write $$O_P=
\left[\begin{array}{cc}
I_{t-2} & 0 \\ 0 & B\\ 
\end{array}\right]\;\;\; \mbox{where}\;\;\; B=\left[\begin{array}{cc}
 & 1\\\;\; b_{ij}\;\; & 0 \\ & \vdots \\ & 0 
\end{array}\right],\;\; b_{ij}=b_{ji}, \;\; 1\leq i,j\leq m+1-t.$$ 
In fact, up to invertible row and column
operations that preserve the almost-symmetry of $B$, we can
assume that $B$ has a symmetric block of maximal size and 
an invertible entry in the last column. Then $$B=\left[\begin{array}{cc}
0 \ldots 0 & 1\\ & 0 \\ b_{ij}  & \vdots \\ & 0 
\end{array}\right],\;\; 2\leq i\leq m+1-t,\; 1\leq j\leq m+1-t,$$
so it contains an almost-symmetric block of size $(m+1-t)\times
(m-t).$ Summarizing, if we let $B'=(b_{ij})_{2\leq i\leq m+1-t,
  1\leq j\leq m+1-t}$ we have that $$O_P=
\left[\begin{array}{cc}
I_{t-1} & 0 \\ 0 & B'\\ 
\end{array}\right].$$ This completes the proof of our claim.

The localization of $I_t(O)$ at the prime ideal $P$ is 
$$I_t(O)_P=I_t(O_P)=I_1(B)$$
thus it is generated by the entries of $B$. Since $B$ is an
almost-symmetric matrix, we have
$$\mu(I_t(O)_P)\leq {m+1-t \choose 2}+(m-t)={m+2-t \choose
  2}-1=c=\hgt I_t(O)_P.$$
Then $I_t(O)$ is locally generated by a regular sequence at all
the minimal associated primes, i.e. $Y$
is generically complete intersection.
\end{enumerate}

\noindent Assume now that $\hgt I_{t-1}(N)=c+1$. Since $I_{t-1}(N)\subseteq I_{t-1}(O)$, then
$$c+1=\hgt I_{t-1}(N)\leq\hgt I_{t-1}(O).$$
Then the two equivalent conditions hold.
\end{proof}

\begin{rmk}
The condition that $I_{t-1}(N)=c+1$ means that $Y$ contains a
symmetric determinantal subscheme $X'$ of codimension 1, whose
defining ideal is
$I_{X'}=I_{t-1}(N).$ Notice that whenever this is the case, $Y$ is
generically complete intersection, hence it is $G_0$. Under this
assumption we have a
concept of generalized divisor on $Y$ (see \cite{ha86a}, \cite{ha94p} and
\cite{ha04u2} about generalized divisors). Then $X'$ is a
generalized divisor on $Y$.  
Theorem~\ref{subm} proves that the existence of such
a subscheme $X'$ of codimension 1 guarantees that $Y$ is locally a complete
intersection. Notice the analogy with standard
determinantal schemes: a standard determinantal scheme $Y$ is good
determinantal if and only if it is locally a complete
intersection, if and only if it contains a standard determinantal
subscheme of codimension 1, whose defining matrix is obtained by deleting a column
from the defining matrix of $Y$.
\end{rmk}

The next example shows that in general the condition that $\hgt
I_{t-1}(N)=c+1$ is stronger than the two equivalent conditions of
Theorem~\ref{subm}.

\begin{ex} Let $R=K[x_0,x_1,x_2,x_3]$, and
let $$O=\left[\begin{array}{ccc}
x_0 & x_1 & x_2 \\
x_1 & x_0 & x_3
\end{array}\right].$$
$I_2(O)$ is a Cohen-Macaulay ideal of height
$2$, $I_1(O)=(x_0,x_1,x_2,x_3)$ has height $4\geq \hgt I_2(O)+1$.
Deleting the last column of $O$ yields the matrix $$N=\left[\begin{array}{cc}
x_0 & x_1 \\
x_1 & x_0 \end{array}\right].$$ $I_1(N)=(x_0,x_1)$ is an ideal of
height $2<2+1.$ Notice that $I_1(N)$ does not change, even if we
perform invertible row and column operations on $O$ that preserve
its almost symmetric structure, before deleting the last column.

Now assume that $O$ is obtained from a homogeneous symmetric
matrix $M$ by deleting the last row. Let $M$ be of the form 
$$M=\left[\begin{array}{ccc}
x_0 & x_1 & x_2 \\
x_1 & x_0 & x_3 \\
x_2 & x_3 & x_2
\end{array}\right].$$
Notice that if we
apply chosen invertible row and column operations that preserve the
symmetry of $M$, we obtain a matrix 
$$M=\left[\begin{array}{ccc}
2x_0+2x_1+2ax_2+2ax_3+a^2x_2 & x_1+x_0+ax_3 & x_2+x_3+ax_2 \\
x_1+x_0+ax_3 & x_0 & x_3 \\
x_2+x_3+ax_2 & x_3 & x_2
\end{array}\right]$$
for any choice of $a\neq 0$.
If we now delete the last row we obtain a new matrix $O'$
$$O'=\left[\begin{array}{ccc}
2x_0+2x_1+2ax_2+2ax_3+a^2x_2 & x_1+x_0+ax_3 & x_2+x_3+ax_2 \\
x_1+x_0+ax_3 & x_0 & x_3
\end{array}\right].$$
Deleting the last column of $O'$ yields the matrix
$$N'=\left[\begin{array}{cc}
2x_0+2x_1+2ax_2+2ax_3+a^2x_2 & x_1+x_0+ax_3 \\
x_1+x_0+ax_3 & x_0
\end{array}\right]$$
and $I_1(N')=(x_0,x_1+ax_3,a(2x_2+ax_2))$. Then $I_1(N')$ has height 3
for any $a\neq 0$. 

Summarizing the example, $I_2(O)$ has maximal height, $\hgt
I_2(O)<\hgt I_1(O)$, $I_1(M)$ has maximal height,  
and $I_1(N)$ does not have maximal
height. However, after applying to $M$ a general transformation that
preserves its symmetry, we obtain $O'$ and $N'$ with the property
that both $I_2(O')$ and $I_1(N')$ have maximal height. 
\end{ex}

The following is the analogous of Theorem~\ref{subm} for
symmetric determinantal schemes. We will not need it in the sequel, 
but we wish to emphasize that a result of this kind holds. 
The proof is very similar to that of the previous theorem, so we
omit it.

\begin{thm}\label{subsd}
Let $X$ be a symmetric determinantal scheme with defining
matrix $M$, $I_X=I_t(M)$. Let $c+1={m-t+2\choose 2}$ 
be the codimension of $X$. The following are equivalent:
\begin{enumerate}
\item $X$ is generically complete intersection.
\item $\hgt I_{t-1}(M)\geq c+2$. 
\end{enumerate}
Let $O$ be the almost-symmetric matrix obtained from $M$ by
  deleting the last row. If $\hgt I_{t-1}(O)\geq c+2$, then 
the two equivalent conditions are verified.
\end{thm}

Finally, we  prove a result in the lines of
the Eisenbud-Evans generalized principal ideal theorem (see
\cite{ei76}) and of its generalization by Bruns (see
\cite{br81a}). The result is not new for an arbitrary matrix $M$,
but the estimate on the height can be 
sharpened in the case of symmetric matrices. We essentially
follow the proof of Theorem~2 in \cite{br81a}.

\begin{thm}\label{ht1}
Let $M$ be a t-homogeneous symmetric matrix of size $m\times m$ with
entries in $R$. Assume that $M$ has no invertible entries. 
Let $O$ be the matrix obtained from $M$ by
deleting the last row, after applying invertible generic row and
column operations to $M$ whcih preserve its symmetry and
t-homogeneity. 
Then $$\hgt I_t(M)/I_t(O)\leq 1.$$
\end{thm}

\begin{proof}
Let $$M=\left[\begin{array}{cccc}
F_{1,1} & \ldots & F_{1,m-1} & F_{1,m} \\
\vdots & & \vdots & \vdots \\
F_{1,m-1} & \ldots & F_{m-1,m-1} & F_{m-1,m} \\
F_{1,m} & \ldots & F_{m-1,m} & F_{m,m}
\end{array}\right]$$ and 
$$O=\left[\begin{array}{cccc}
F_{1,1} & \ldots & F_{1,m-1} & F_{1,m} \\
\vdots & & \vdots & \vdots \\
F_{1,m-1} & \ldots & F_{m-1,m-1} & F_{m-1,m}
\end{array}\right].$$
Consider the matrix $$L=\left[\begin{array}{ccccc}
F_{1,1} & \ldots & F_{1,m-1} & F_{1,m} & 0 \\
\vdots & & \vdots & \vdots & \vdots \\
F_{1,m-1} & \ldots & F_{m-1,m-1} & F_{m-1,m} & 0 \\
F_{1,m} & \ldots & F_{m-1,m} & F_{m,m} & -1
\end{array}\right].$$ Notice that $L$ and $M$ are related 
in the same way as $\varphi$ and $\varphi'$ in the proof of 
Theorem 2 of \cite{br81a}. We regard the matrices over the ring
$S=R/I_t(O)$. $L$ defines a morphism $\psi:S^m\lra S^{m+1}$
where the images of a basis of $S^m$ are given by the rows of
$L$. Similarly, $M$ defines a morphism $\varphi:S^m\lra S^m$.
Let $\M:=Coker \psi$ and $\M':=Coker \varphi$. Then $\M'\cong
\M/e_{m+1}$, where $e_1,\ldots,e_{m+1}$ denote the elements of the
standard basis of $S^{m+1}$. For an $x\in\M$, we define 
$$\M^*(x):=\{ f(x)\; : \; f\in Hom_S(\M,S) \}.$$
As shown in Theorem~2 of~\cite{br81a}, we have that 
$$I_t(M)/I_t(O)\subseteq \M^*(e_{m+1}).$$
Following the proof of Theorem 2 of \cite{br81a}, one can show
that $\hgt \M^*(e_{m+1})\leq m-t+1$. However we claim that in
the special case of a symmetric matrix $M$, one has the sharper 
bound $\hgt \M^*(e_{m+1})\leq 1$.

Since $e_{m+1}=F_{1,m}e_1+\ldots+F_{m,m}e_m$, we have 
$$\M^*(e_{m+1})=\{  f(F_{1,m}e_1+\ldots+F_{m,m}e_m)\; : \; f\in
Hom_S(\M,S) \}.$$ For each $f\in Hom_S(\M,S)$, we have
$$f(F_{1,m}e_1+\ldots+F_{m,m}e_m)=F_{1,m}f(e_1)+\ldots+F_{m,m}f(e_m).$$
Let $f_i:=f(e_i)\in S$.
Assume by contradiction that $\hgt \M^*(e_{m+1})\geq 2$. Then we
can find $f,g\in Hom_S(\M,S)$ such that $f(e_{m+1})$ and
$g(e_{m+1})$ form a regular sequence in $S$. For all
$i=1,\ldots,m-1$ we have 
$$f(\sum_{j=1}^m F_{i,j}e_j)=\sum_{j=1}^m F_{i,j}f_j=0$$
therefore $$f_mF_{i,m}=-\sum_{j=1}^{m-1}F_{i,j}f_j.$$
Analogously, $$g_mF_{i,m}=-\sum_{j=1}^{m-1}F_{i,j}g_j.$$
Hence $$g_mf(e_{m+1})=\sum_{i=1}^m
F_{i,m}f_ig_m=f_mg_mF_{m,m}-\sum_{\begin{array}{c} i=1,\ldots,m-1
    \\ j=1,\ldots,m-1\end{array}} f_ig_jF_{i,j}$$ and since $M$
is symmetric
$$g_mf(e_{m+1})=f_mg(e_{m+1}).$$ By assumption $f(e_{m+1})$ and
$g(e_{m+1})$ form a regular sequence in $S$, so $g_m=hg(e_{m+1})$
and $f_m=hf(e_{m+1})$ for some $h\in S$. Assume that
$f_m,g_m\neq 0$. So $h\neq 0$ and \begin{equation}\label{degr}
f_m=hf(e_{m+1})=\sum_{i=1}^m
hF_{i,m}f_i=\sum_{i\in I}
hF_{i,m}f_i\end{equation}
where $I\subseteq \{1,\ldots,m\}$ is the set of indexes
of the summands that effectively contribute to the sum. In other
words, if $f=(f_1,\ldots,f_m)$ and $\sum_{i\not\in I}
F_{i,m}f_i=0$ then $f(e_{m+1})=\phi(e_{m+1})$ where $\phi=(\phi_1,\ldots,\phi_m)$
with $\phi_i=f_i$ if $i\in I$ and $\phi_i=0$ if $i\not\in I$. So
we can replace $f$ with $\phi$. This proves that we can assume
without loss of generality that since $f_m\neq 0$, then $m\in I$.
Hence the term $F_{m,m}f_m$ appears in the sum (\ref{degr}).
Comparing degrees in (\ref{degr}), we get 
$$deg(f_m)\geq deg(h)+deg(F_{m,m})+deg(f_m)\geq
deg(F_{m,m})+deg(f_m),$$ hence $deg(F_{m,m})\leq 0.$ By the
assumption that $M$ has no invertible entries, $F_{m,m}=0$. But
this is a contradiction: since we are allowing generic invertible
row and column operations that preserve the symmetry of $M$ we
can always assume that $F_{m,m}\neq 0$ unless all the entries in
the last row and column of $M$ are zero. However, in that case 
$I_t(O)=I_t(M)$ and the thesis is trivially verified.
Notice if $K$ has characteristic $2$ and $F_{m,m}=0$, any row and
column operations that preserve the symmetry of $M$ would also
preserve the property that $F_{m,m}=0$.

We still need to analyze the case when $f_m=g_m=0$. We have
$$f(\sum_{j=1}^m F_{i,j}e_j)=\sum_{j=1}^{m-1} F_{i,j}f_j=0,$$
therefore $$f_{m-1}F_{i,m-1}=-\sum_{j=1}^{m-2}F_{i,j}f_j.$$
Analogously, $$g_{m-1}F_{i,m-1}=-\sum_{j=1}^{m-2}F_{i,j}g_j.$$
Hence, proceeding as in the previous case, 
$$f_{m-1}g(e_{m+1})=g_{m-1}f(e_{m+1})$$ and either
$f_{m-1},g_{m-1}\neq 0$ or $f_{m-1}=g_{m-1}=0$. In the first
case, we can conclude as above, in the second case we obtain 
$$f_{m-2}g(e_{m+1})=g_{m-2}f(e_{m+1}).$$ We can keep iterating
this reasoning until either $f_i,g_i\neq 0$
for some $i$, or $f=g=0$. In both cases we get a contradiction. 
\end{proof}

\section{Biliaison of symmetric determinantal schemes}

In this section we prove that symmetric determinantal schemes
are in the same G-biliaison class of a complete intersection
of the same codimension.
We start by proving that any symmetric
determinantal scheme is a divisor on an almost-symmetric
determinantal scheme.
This result will be used in the proof
of Theorem~\ref{bilsym}.

\begin{prop}\label{gsd}
Let $X$ be a symmetric determinantal scheme of codimension $c+1$. Assume that 
the ideal of $X$ is generated by the $t\times t$ minors of a
t-homogeneous matrix $M$, $I_X=I_t(M)$.
Let $O$ be the matrix obtained from $M$ by deleting the last row
(after performing generic invertible row and column operations
that preserve the symmetry and homogeneity of $M$). 
Let $N$ be the matrix obtained from $O$ by deleting the last column.
Then:\begin{itemize} 
\item $N$ is a t-homogeneous symmetric matrix. It defines a symmetric
determinantal scheme $X'$ of codimension $c+1$, with
$I_{X'}=I_{t-1}(N)$.
\item $O$ is a t-homogeneous almost-symmetric matrix. It defines an almost-symmetric
determinantal scheme $Y$ of codimension $c$, with
$I_{Y}=I_{t}(O)$. $Y$ is an arithmetically Cohen-Macaulay, 
generically complete intersection scheme. 
\end{itemize}
\end{prop}

\begin{proof}
Let $N$ be the matrix obtained from $M$ by deleting a row and
the corresponding column, after performing generic invertible row
and column operations on $M$ that preserve its symmetry and
t-homogeneity. Then $N$
is symmetric, and $\hgt I_{t-1}(N)\leq
{m-1-(t-1)+2\choose 2}=c+1$ by Remark~\ref{codim}. 
Let $O$ be the matrix obtained from $M$ by deleting the last row,
$O\supseteq N$.
By Theorem~\ref{ht1}, $$\hgt I_t(M)/I_t(O)\leq 1,$$
then $$\hgt I_t(O)\geq \hgt I_t(M)-1=c.$$ It follows that $\hgt
I_t(O)=c$ and $I_t(O)$ is Cohen-Macaulay by Theorem~\ref{gsd1}.
Therefore $I_Y=H^0_*(I_t(O))=I_t(O),$ and the scheme $Y$ is
arithmetically Cohen-Macaulay of codimension $c$. 
Moreover, $$\hgt I_{t-1}(O)\geq \hgt I_t(M)=c+1>c=\hgt I_t(O),$$
so $Y$ is generically complete intersection by Theorem~\ref{subm}.

Consider the ideals $I_t(M)/I_t(O)$ and $I_{t-1}(N)/I_t(O)$ contained
in $R/I_t(O)$. $R/I_t(O)$ is a Cohen-Macaulay ring, and
$I_t(M)/I_t(O)$ is an ideal of height 1. 
Let $R'$ be the ring of total quotients of $R/I_t(O)$, i.e. the
localization of $R/I_t(O)$ at the set $S$ consisting of all its 
nonzerodivisors. $I_t(M)/I_t(O)$
and $I_{t-1}(N)/I_t(O)$ are isomorphic as submodules of $R'$, the
isomorphism being given by multiplication by 
$M_{i_1,\ldots,i_{t-1};j_1,\ldots,j_{t-1}}/M_{i_1,\ldots,i_{t-1},m;j_1,\ldots,j_{t-1},m}$ 
for any choice of $1\leq i_1<\ldots<i_{t-1}<m$ and $1\leq j_1<\ldots<j_{t-1}<m $. 
Here $M_{i_1,\ldots,i_{t-1},m;j_1,\ldots,j_{t-1},m}$ denotes the minor
of $M$ that corresponds to rows $i_1,\ldots,i_{t-1},m$ and columns
$j_1,\ldots,j_{t-1},m$, and $M_{i_1,\ldots,i_{t-1};j_1,\ldots,j_{t-1}}$ denotes the minor
of $M$ that corresponds to rows $i_1,\ldots,i_{t-1}$ and columns
$j_1,\ldots,j_{t-1}$. The inverse is given by multiplication by 
$M_{i_1,\ldots,i_{t-1},m;j_1,\ldots,j_{t-1},m}/M_{i_1,\ldots,i_{t-1};j_1,\ldots,j_{t-1}}.$
See the proof of Theorem~\ref{bilsym} for more details about the isomorphism.
Then $\hgt I_{t-1}(N)/I_t(O)=1$, hence $\hgt I_{t-1}(N)=c+1$. 
So $I_{t-1}(N)$ defines an
arithmetically Cohen-Macaulay, symmetric determinantal scheme
$X'$ of codimension $c+1$, and $I_{t-1}(N)=I_{X'}$.
\end{proof}

\begin{rmk}
If $X$ is a symmetric determinantal scheme of codimension $c+1$ with defining
matrix $M$, $I_X=I_t(M)$, we prove that:
\begin{itemize}
\item $\hgt I_{t-1}(O)\geq c+1$, where $O$ is almost-symmetric and obtained 
from $M$ by deleting the last row.
\item $\hgt I_{t-1}(N)= c+1$, where $N$ is symmetric and obtained from $M$ by
  deleting the last row and column.
\end{itemize}
Notice that by Theorem~\ref{subm} this implies that the scheme
$Y$ defined by $I_t(O)$ is generically complete
intersection. However, it does not imply the same result for $X$,
nor for the scheme $X'$ defined by $I_{t-1}(N)$. 
\end{rmk}

We are now ready to prove the main result of this paper.

\begin{thm}\label{bilsym}
Any symmetric determinantal scheme in $\PP^n$ can be obtained from a linear
variety by a finite sequence of ascending elementary G-biliaisons.
\end{thm}

\begin{proof}
Let $X\subset\PP^n$ be a symmetric determinantal scheme. We
follow the notation of Definition~\ref{sd}. Let $M=(F_{ij})$ be
the matrix whose minors of size $t\times t$ define $X$. From the definition
$F_{ij}=F_{ji}$ for all $i,j$, and the matrix is t-homogeneous. 
Let $c+1$ be the codimension of
$X$, $c:={m-t+2 \choose 2}-1$. 
If $t=1$ or $t=m$, then $X$ is a complete intersection, therefore we
can perform a finite sequence of descending elementary
CI-biliaisons to a linear variety. Therefore, we concentrate on
the case when $2\leq t<m$.

Let $O$ be the matrix obtained from $M$ by deleting the last
row, after performing generic invertible row and column
operations that preserve the symmetry of $M$. 
$O$ is a t-homogeneous matrix of size $(m-1)\times m$.
Let $Y$ be the scheme whose saturated ideal is generated by the
$t\times t$ minors of $O$. By Proposition~\ref{gsd}, $Y$ is an arithmetically 
Cohen-Macaulay scheme of codimension $c$. Notice that $Y$ is
standard determinantal exactly when $t=m-1$. For our purpose, it is important 
to observe that $Y$ is generically complete intersection. In particular, it 
satisfies the property $G_0$. Therefore, a biliaison on 
$Y$ is a G-biliaison, hence also an even G-liaison. This was proved in 
\cite{kl01} for $Y$ satisfying property $G_1$ and extended in \cite{ha04u2}
to $Y$ satisfying property $G_0$.

Let $N$ be the matrix obtained from $M$ by deleting the last row
and column. $N$ is a t-homogeneous symmetric matrix of size 
$(m-1)\times(m-1)$.
Let $X'$ be the scheme cut out by the
$(t-1)\times(t-1)$ minors of $N$. 
Both $X$ and $X'$ are contained in $Y$.
We denote by $H$ a hyperplane section divisor on $Y$.
We are going to show that $$X\sim X'+aH \;\;\;\mbox{for some $a>0$},$$ 
where $\sim$ denotes
linear equivalence of divisors. This will prove that $X$ is
obtained by an elementary biliaison from $X'$.
Continuing in this manner, after $t-1$ biliaisons we
reduce to the case $t=1$, when the scheme $X$ is a complete
intersection. Then we can perform descending CI-biliaisons to a
linear variety.

Let $\I_{X|Y}$, $\I_{X'|Y}$  be the ideal sheafs on $Y$ of $X$ and
$X'$. We then need to show that 
\begin{equation}\label{isosh} 
\I_{X|Y}\cong \I_{X'|Y}(-a)\;\;\;\mbox{for some $a>0$}.
\end{equation}

A system of generators of $I_{X|Y}=H^0_*(\I_{X|Y})=I_t(M)/I_Y$ is
given by the images in the coordinate ring of $Y$ of the $t\times
t$ minors of $M$
$$I_{X|Y}=(M_{i_1,\ldots,i_t;j_1,\ldots,j_t} \; | \; 1\leq
i_1<i_2<\ldots<i_t\leq m, 1\leq
j_1<j_2<\ldots<j_t\leq m).$$
Here $M_{i_1,\ldots,i_t;j_1,\ldots,j_t}$ denotes the image of
the determinant of the submatrix of $M$ consisting of rows
$i_1,\ldots,i_t$ and columns $j_1,\ldots,j_t$ in
the coordinate ring of $Y$. The saturated ideal of $Y$ is
minimally generated by the minors of size $t\times t$ of $M$ that do not
involve the last row. 
Notice that this minimal system of generators of
$I_Y$ can be completed to a minimal system of generators of $I_X$
by adding all the minors of size $t\times t$ of $M$ that involve
both the last row and the last column.
Therefore, a minimal system of generators of
$I_{X|Y}$ is given by 
$$I_{X|Y}=(M_{i_1,\ldots,i_{t-1},m;j_1,\ldots,j_{t-1},m} \; | \; 1\leq
i_1<\ldots<i_{t-1}\leq m-1, 1\leq j_1<\ldots<j_{t-1}\leq m-1).$$

A minimal system of generators of $I_{X'|Y}=H^0_*(\I_{X'|Y})=
I_{t-1}(N)/I_Y$ is
given by the images in the coordinate ring of $Y$ of the $(t-1)\times
(t-1)$ minors of $N$ 
$$I_{X'|Y}=(M_{i_1,\ldots,i_{t-1};j_1,\ldots,j_{t-1}} \; | \; 1\leq
i_1<\ldots<i_{t-1}\leq m-1, 1\leq
j_1<\ldots<j_{t-1}\leq m-1).$$
Here $M_{i_1,\ldots,i_{t-1};j_1,\ldots,j_{t-1}}$ denotes the image of
the determinant of the submatrix of $M$ consisting of rows
$i_1,\ldots,i_{t-1}$ and columns $j_1,\ldots,j_{t-1}$ in
the coordinate ring of $Y$.

In order to prove the isomorphism (\ref{isosh}), it suffices to
check that the quotients 
\begin{equation}\label{eqk}
\frac{M_{i_1,\ldots,i_{t-1},m;j_1,\ldots,j_{t-1},m}}
{M_{i_1,\ldots,i_{t-1};j_1,\ldots,j_{t-1}}}
\end{equation}
are all equal as elements of $H^0(\K_Y(a))$, where $\K_Y$ is the
sheaf of total quotient rings of $Y$. This also gives us an easy
way to compute the value of $a$ as the difference 
$deg(M_{i_1,\ldots,i_{t-1},m;j_1,\ldots,j_{t-1},m})-
deg(M_{i_1,\ldots,i_{t-1};j_1,\ldots,j_{t-1}})=deg(F_{m,m}).$

Equality (\ref{eqk}) is readily verified, once we show that
$$M_{i_1,\ldots,i_{t-1},m;j_1,\ldots,j_{t-1},m}\cdot
M_{k_1,\ldots,k_{t-1};l_1,\ldots,l_{t-1}}-
M_{k_1,\ldots,k_{t-1},m;l_1,\ldots,l_{t-1},m}\cdot
M_{i_1,\ldots,i_{t-1};j_1,\ldots,j_{t-1}}\in
I_Y.$$ The proof is then completed by the following lemmas.
\end{proof}

Since we could not find an adequate reference in the literature, we need to 
prove the following two lemmas about the minors of a matrix.

\begin{lemma}\label{matx}
Let $M$ be a matrix of size $m\times m$ and let
$M_{i_1,\ldots,i_a;j_1,\ldots,j_a}$ denote the minor of the
submatrix of $M$ consisting of rows $i_1,\ldots,i_a$ and columns
$j_1,\ldots,j_a$. Let $I$ be the ideal generated by the minors of
$M$ of size $(a+1)\times (a+1)$.
Then
$$M_{i_1,\ldots,i_a;j_1,\ldots,j_a}\cdot
M_{k_1,\ldots,k_a;l_1,\ldots,l_a}-
M_{k_1,\ldots,k_a;j_1,\ldots,j_a}\cdot
M_{i_1,\ldots,i_a;l_1,\ldots,l_a}\in I.$$
\end{lemma}

\begin{proof}
We start by proving the thesis when $i_b=k_
b$ for $b=1,\ldots a-1$. So we want to show that
\begin{equation}\label{intermediate}
M_{i_1,\ldots,i_a;j_1,\ldots,j_a}\cdot
M_{i_1,\ldots,i_{a-1},k_a;l_1,\ldots,l_a}-
M_{i_1,\ldots,i_{a-1},k_a;j_1,\ldots,j_a}\cdot
M_{i_1,\ldots,i_a;l_1,\ldots,l_a}\in I.\end{equation}
This is essentially an application
of Sylvester's identity:
$$M_{i_1,\ldots,i_a;j_1,\ldots,j_a}\cdot
M_{i_1,\ldots,i_{a-1},k_a;j_1,\ldots,j_{a-1},l_a}-
M_{i_1,\ldots,i_{a-1},k_a;j_1,\ldots,j_a}\cdot
M_{i_1,\ldots,i_a;j_1,\ldots,j_{a-1},l_a}=$$
$$M_{i_1,\ldots,i_{a-1};j_1,\ldots,j_{a-1}}\cdot
M_{i_1,\ldots,i_a,k_a;j_1,\ldots,j_a,l_a}$$
For our purpose, we only need that the difference belongs to $I$.
See \cite{mu60}, pg. 33, for a general statement and a proof of
Sylvester's identity.
$$M_{i_1,\ldots,i_a;j_1,\ldots,j_a}\cdot
M_{i_1,\ldots,i_{a-1},k_a;l_1,\ldots,l_a}-
M_{i_1,\ldots,i_{a-1},k_a;j_1,\ldots,j_a}\cdot
M_{i_1,\ldots,i_a;l_1,\ldots,l_a}=$$ $$M_{i_1,\ldots,i_a;j_1,\ldots,j_a}\cdot
M_{i_1,\ldots,i_{a-1},k_a;l_1,\ldots,l_a}-\sum_{b=1}^{a-1}
M_{i_1,\ldots,i_a;j_1,\ldots,j_b,l_{b+1},
\ldots,l_a}\cdot M_{i_1,\ldots,i_{a-1},k_a;l_1,\ldots,l_b,j_{b+1},\ldots,j_a}
+$$ $$\sum_{b=1}^{a-1} M_{i_1,\ldots,i_a;j_1,\ldots,j_b,l_{b+1},
\ldots,l_a}\cdot M_{i_1,\ldots,i_{a-1},k_a;l_1,\ldots,l_b,j_{b+1},\ldots,j_a}
-M_{i_1,\ldots,i_{a-1},k_a;j_1,\ldots,j_a}\cdot
M_{i_1,\ldots,i_a;l_1,\ldots,l_a}=$$
$$=\sum_{b=1}^{a} (M_{i_1,\ldots,i_a;j_1,\ldots,j_b,l_{b+1},
\ldots,l_a}\cdot
M_{i_1,\ldots,i_{a-1},k_a;l_1,\ldots,l_b,j_{b+1},\ldots,j_a}-$$
$$M_{i_1,\ldots,i_a;j_1,\ldots,j_{b-1},l_b,
\ldots,l_a}\cdot
M_{i_1,\ldots,i_{a-1},k_a;l_1,\ldots,l_{b-1},j_b,\ldots,j_a}).$$
Each summand is of the form $$M_{i_1,\ldots,i_a;j_1,\ldots,j_b,l_{b+1},
\ldots,l_a}\cdot
M_{i_1,\ldots,i_{a-1},k_a;l_1,\ldots,l_b,j_{b+1},\ldots,j_a}-$$
$$M_{i_1,\ldots,i_a;j_1,\ldots,j_{b-1},l_b,\ldots,l_a}\cdot
M_{i_1,\ldots,i_{a-1},k_a;l_1,\ldots,l_{b-1},j_b,\ldots,j_a},$$
for $b=1,\ldots,a$.
In particular, all the minors in the expression have all the rows and
the columns in common, except possibly for one.
So Sylvester's identity applies, and the thesis follows.

Let us now prove the thesis in full generality. We are going to
use (\ref{intermediate}), and we will proceed in an analogous
manner to the proof above.
We want to show that 
$$M_{i_1,\ldots,i_a;j_1,\ldots,j_a}\cdot
M_{k_1,\ldots,k_a;l_1,\ldots,l_a}-
M_{k_1,\ldots,k_a;j_1,\ldots,j_a}\cdot
M_{i_1,\ldots,i_a;l_1,\ldots,l_a}\in I.$$
Rewrite the difference as $$M_{i_1,\ldots,i_a;j_1,\ldots,j_a}\cdot
M_{k_1,\ldots,k_a;l_1,\ldots,l_a}-
M_{k_1,\ldots,k_a;j_1,\ldots,j_a}\cdot
M_{i_1,\ldots,i_a;l_1,\ldots,l_a}=$$
$$M_{i_1,\ldots,i_a;j_1,\ldots,j_a}\cdot
M_{k_1,\ldots,k_a;l_1,\ldots,l_a}-
\sum^{b=1}_{a-1} M_{i_1,\ldots,i_b,k_{b+1},\ldots,k_a;j_1,\ldots,j_a}\cdot
M_{k_1,\ldots,k_b,i_{b+1},\ldots,i_a;l_1,\ldots,l_a}+$$
$$\sum_{b=1}^{a-1} M_{i_1,\ldots,i_b,k_{b+1},\ldots,k_a;j_1,\ldots,j_a}\cdot
M_{k_1,\ldots,k_b,i_{b+1},\ldots,i_a;l_1,\ldots,l_a}
-M_{k_1,\ldots,k_a;j_1,\ldots,j_a}\cdot
M_{i_1,\ldots,i_a;l_1,\ldots,l_a}=$$
$$\sum_{b=1}^a M_{i_1,\ldots,i_b,k_{b+1},\ldots,k_a;j_1,\ldots,j_a}\cdot
M_{k_1,\ldots,k_b,i_{b+1},\ldots,i_a;l_1,\ldots,l_a}-$$
$$M_{i_1,\ldots,i_{b-1},k_b,
\ldots,k_a;j_1,\ldots,j_a}\cdot
M_{k_1,\ldots,k_{b-1},i_b,\ldots,i_a;l_1,\ldots,l_a}.$$
Each summand is of the form $$M_{i_1,\ldots,i_b,k_{b+1},\ldots,k_a;j_1,\ldots,j_a}\cdot
M_{k_1,\ldots,k_b,i_{b+1},\ldots,i_a;l_1,\ldots,l_a}-
M_{i_1,\ldots,i_{b-1},k_b,
\ldots,k_a;j_1,\ldots,j_a}\cdot
M_{k_1,\ldots,k_{b-1},i_b,\ldots,i_a;l_1,\ldots,l_a},$$
hence by (\ref{intermediate}) it belongs to $I$.
This concludes the proof.
\end{proof}

\begin{rmk}
Notice that following the steps of the proof, one can write down
explicitly $M_{i_1,\ldots,i_a;j_1,\ldots,j_a}\cdot
M_{k_1,\ldots,k_a;l_1,\ldots,l_a}-
M_{k_1,\ldots,k_a;j_1,\ldots,j_a}\cdot
M_{i_1,\ldots,i_a;l_1,\ldots,l_a}$ as a combination of the minors
of size $(a+1)\times (a+1)$ of the matrix $M$. This is not
relevant for our purposes.
\end{rmk}

The following lemma concludes the proof of Theorem~\ref{bilsym}.

\begin{lemma}
Let $M$ be a matrix of size $m\times m$ and let
$M_{i_1,\ldots,i_t;j_1,\ldots,j_t}$ denote the minor of the
submatrix of $M$ consisting of rows $i_1,\ldots,i_t$ and columns
$j_1,\ldots,j_t$. Let $I_Y$ be the ideal generated by the determinants
of the submatrices of $M$ of size $t\times t$ of $M$, that do not
contain the last row. Then
$$M_{i_1,\ldots,i_{t-1},m;j_1,\ldots,j_{t-1},m}\cdot
M_{k_1,\ldots,k_{t-1};l_1,\ldots,l_{t-1}}-
M_{k_1,\ldots,k_{t-1},m;l_1,\ldots,l_{t-1},m}\cdot
M_{i_1,\ldots,i_{t-1};j_1,\ldots,j_{t-1}}\in I_Y.$$
\end{lemma}

\begin{proof}
It is enough to prove that the statement holds for a matrix of indeterminates
$M=(x_{ij})$, $1\leq i\leq j\leq m$.
From Lemma~\ref{matx}, we have that 
\begin{equation}\label{eqsymm}
M_{i_1,\ldots,i_{t-1},m;j_1,\ldots,j_{t-1},m}\cdot
M_{k_1,\ldots,k_{t-1},m;l_1,\ldots,l_{t-1},m}-
M_{k_1,\ldots,k_{t-1},m;j_1,\ldots,j_{t-1},m}\cdot
M_{i_1,\ldots,i_{t-1},m;l_1,\ldots,l_{t-1},m}\in I\subseteq I_Y.
\end{equation}
Here $I$ is the ideal generated by the minors of $M$ of size $t+1$.
Clearly, $I\subseteq I_Y$.

Expanding the determinant $M_{i_1,\ldots,i_{t-1},m;j_1,\ldots,j_{t-1},m}$
about column $m$ we obtain
$$M_{i_1,\ldots,i_{t-1},m;j_1,\ldots,j_{t-1},m}=\sum_{h=1}^t (-1)^h x_{i_h,m}
M_{i_1,\ldots,i_{h-1},i_{h+1},\ldots,i_{t-1},m;j_1,\ldots,j_{t-1}}.$$
We are  adopt the convention that $i_t=j_t=k_t=l_t=m$.

Substituting this expression in the equation (\ref{eqsymm}), we get 
$$M_{i_1,\ldots,i_{t-1},m;j_1,\ldots,j_{t-1},m}\cdot
M_{k_1,\ldots,k_{t-1},m;l_1,\ldots,l_{t-1},m}-
M_{k_1,\ldots,k_{t-1},m;l_1,\ldots,l_{t-1},m}\cdot
M_{i_1,\ldots,i_{t-1},m;j_1,\ldots,j_{t-1},m}=$$
$$\sum_{h=1}^t (-1)^h x_{i_h,m}
M_{i_1,\ldots,i_{h-1},i_{h+1},\ldots,i_{t-1},m;j_1,\ldots,j_{t-1}}
M_{k_1,\ldots,k_{t-1},m;l_1,\ldots,l_{t-1},m}-$$
$$\sum_{h=1}^t (-1)^h x_{k_h,m}
M_{k_1,\ldots,k_{h-1},k_{h+1},\ldots,k_{t-1},m;l_1,\ldots,l_{t-1}}
M_{i_1,\ldots,i_{t-1},m;j_1,\ldots,j_{t-1},m}\in I_S$$

The coefficient of $x_{m,m}$ in (\ref{eqsymm}) is then
$$(-1)^t
M_{i_1,\ldots,i_{t-1};j_1,\ldots,j_{t-1}}
M_{k_1,\ldots,k_{t-1},m;l_1,\ldots,l_{t-1},m} - (-1)^t
M_{k_1,\ldots,k_{t-1};l_1,\ldots,l_{t-1}}
M_{i_1,\ldots,i_{t-1},m;j_1,\ldots,j_{t-1},m}.$$
Since (\ref{eqsymm}) is an equation of the form 
$\alpha x_{m,m}+\beta=0$ mod. $I_Y$, and since this equality has 
to hold for any $x_{m,m}$, we deduce that $\alpha=0$
mod. $I_Y$. Equivalently, $$M_{i_1,\ldots,i_{t-1};j_1,\ldots,j_{t-1}}
M_{k_1,\ldots,k_{t-1},m;l_1,\ldots,l_{t-1},m} - 
M_{k_1,\ldots,k_{t-1};l_1,\ldots,l_{t-1}}
M_{i_1,\ldots,i_{t-1},m;j_1,\ldots,j_{t-1},m}\in I_Y,$$
that is what we wanted to prove.
\end{proof}

We want to emphasize a consequence of the proof of Theorem~\ref{bilsym}.

\begin{cor}\label{glicci}
Every symmetric determinantal scheme $X$ can be G-bilinked in $t-1$ steps to a 
complete intersection, where $t$ is the size of the minors defining $X$.
In particular, every symmetric determinantal scheme is glicci.
\end{cor}

We end with an example that shows that the proof of
Theorem~\ref{bilsym} does not extend to a field $K$ of characteristic
2. The matrix $M$ below was brought to our attention by W. Bruns
as an example of a matrix such that $\hgt I_2(M)=3$ is maximal, but
$I_3(M)=0$ if the field $K$ has characteristic $2$.

\begin{ex}
Let $K$ be an algebraically closed field of characteristic $2$, let
$X\subseteq\PP^3$ be the fat point with $I_X=(x,y,z)^2\subseteq
R=K[x,y,z,w]$. $X$ is symmetric determinantal, since $I_X=I_2(M)$
where 
$$M=\left[\begin{array}{ccc} 0 & x & y \\
x & 0 & z \\ y & z & 0 \end{array}\right]$$ and $\hgt
I_2(M)=3={3-2+2\choose 2}$. Any invertible row and column operations
that preserve the symmetry and homogeneity of $M$ yield a matrix $M'$
of the form $$M'=\left[\begin{array}{ccc} 0 & L_1 & L_2 \\
L_1 & 0 & L_3 \\ L_2 & L_3 & 0 \end{array}\right]$$ where
$L_1,L_2,L_3$ are linearly independent linear forms. So the matrix $O$
has the form $$O=\left[\begin{array}{ccc} 0 & L_1 & L_2 \\
L_1 & 0 & L_3 \end{array}\right]$$ and $\hgt I_2(O)=1$. Moreover,
deleting the last column of $O$ we obtain $$N=\left[\begin{array}{cc} 0 & L_1 \\
L_1 & 0 \end{array}\right]$$ and $\hgt I_2(N)=1$. Therefore, for $X$
we cannot produce a G-biliaison as discussed in Theorem~\ref{bilsym}.
Notice however that $X$ is
standard determinantal, corresponding to the matrix 
$$\left[\begin{array}{cccc} x & y & z & 0 \\
0 & x & y & z\end{array}\right].$$ Hence $X$ can be
obtained from the point $[0:0:0:1]$ via a G-biliaison.
\end{ex}

\end{document}